\documentclass[12pt]{amsart}

\title[]{Existence of  divergent Birkhoff normal forms of Hamiltonian functions}

\thanks{Supported in part by NSF grant DMS-0305474}
\author[]{Xianghong Gong}
\address{Department of Mathematics,
University of Wisconsin, Madison, WI 53706}
\email{gong@math.wisc.edu}

\subjclass[2000]{Primary 37J40}

\newtheorem{thm}{Theorem}[section]
\newtheorem{cor}[thm]{Corollary}
\newtheorem{prop}[thm]{Proposition}
\newtheorem{lemma}[thm]{Lemma}

\theoremstyle{definition}

\newcommand{\cc}{{\bf C}}

\newcommand{\rr}{{\bf R}}

\newcommand{\ov}{\overline}

\newcommand{\bg}{\begin{gather}}
\newcommand{\bgn}{\begin{gather*}}
\newcommand{\egn}{\end{gather*}}
\newcommand{\be}[1]{\begin{equation}\label{#1}}
\newcommand{\ben}{\begin{equation*}}
\newcommand{\beq}{\begin{eqnarray}}
\newcommand{\eeq}{\end{eqnarray}}
\newcommand{\bc}[1]{\begin{cor}\label{cor:#1}}
\newcommand{\ec}{\end{cor}}
\newcommand{\bt}[1]{\begin{thm}\label{thm:#1}}
\newcommand{\et}{\end{thm}}
\newcommand{\bl}[1]{\begin{lemma}\label{lemma:#1}}
\newcommand{\el}{\end{lemma}}
\newcommand{\bd}[1]{\begin{define}\label{def:#1}}
\newcommand{\ed}{\end{define}}
\newcommand{\br}[1]{\begin{remark}\label{def:#1}}
\newcommand{\er}{\end{remark}}
\newcommand{\bp}[1]{\begin{prop}\label{prop:#1}}
\newcommand{\ep}{\end{prop}}
\newcommand{\bs}{\begin{proof}[Solution]}
\newcommand{\es}{\end{proof}}

\newcommand{\re}[1]{(\ref{#1})}
\newcommand{\rea}[1]{$(\ref{#1})$}
\newcommand{\rl}[1]{Lemma~$\ref{lemma:#1}$}

\newcommand{\rt}[1]{Theorem~$\ref{thm:#1}$}

\newcounter{pp}
\newcommand{\bpp}{\begin{list}{$\hspace{-1em}\arabic{pp}.$}{\usecounter{pp}}}
\newcommand{\epp}{\end{list}}

\newcounter{ppp}
\newcommand{\bppp}{\begin{list}{$\hspace{-1em}(\alph{ppp})$}{\usecounter{ppp}}}
\newcommand{\eppp}{\end{list}}

\newcommand{\Pezeth}{\bibitem{Pezeth} R. P\'erez-Marco,
{\em Convergence or generic divergence of the Birkhoff normal form},
Ann. of Math. (2) {\bf 157}(2003), no. 2, 557-574.
}

\newcommand{\Elnize}{\bibitem{Elnize}
L.H. Eliasson,
{\em Normal forms for Hamiltonian systems with Poisson commuting integrals--elliptic case},
Comment. Math. Helv.  {\bf 65}(1990),  no. 1, 4--35.
}

\newcommand{\Gizeon}{\bibitem{Gizeon} A. Giorgilli,
 {\em Unstable equilibria of Hamiltonian systems},
  Discrete Contin. Dynam. Systems  {\bf 7}(2001),  no. 4, 855--871
  }

\newcommand{\zezemrl}{\bibitem{zezemrl} X. Gong,
{\em Conformal maps, monodromy transformations, and non-reversible Hamiltonian
systems},  Math. Res. Lett. {\bf 7}(2000), no.~4, 471-476.}

\newcommand{\zeonaif}{\bibitem{zeonaif}  X. Gong,
{\em Levi-flat invariant sets of holomorphic symplectic mappings},
Ann. Inst. Fourier(Grenoble) {\bf 51}(2001), no.~1, 151-208. }

\newcommand{\Rusifo}{\bibitem{Rusifo} H. R\"{u}ssmann,
{\em
\"Uber das Verhalten analytischer Hamiltonscher
Differentialgleichungen in der N\"ahe einer Gleichgewichtsl\"osung},
Math. Ann. {\bf 154}(1964), 285--300.
}

\newcommand{\Rusise}{\bibitem{Rusise}H. R\"{u}ssmann,
{\em \"{U}ber die Normalform analytischer Hamiltonscher
Differentialgleichungen in der N\"ahe einer Gleichgewichtsl\"{o}sung},
Math. Ann. {\bf 169}(1967), 55-72.
}

\newcommand{\Veseei}{\bibitem{Veseei} J. Vey,
{\em Sur certains syst\`{e}mes dynamiques s\'{e}perables},
Amer. J. Math. {\bf 100} (1978), 591-614.
}

\newcommand{\Sifoon}{\bibitem{Sifoon} C.L. Siegel,
{\em On integrals of canonical systems},
Ann. Math. {\bf 42}(1941), 806-822.
}

\newcommand{\Iteini}{\bibitem{Iteini} H. Ito,
{\em Convergence of Birkhoff normal forms for integrable systems},
Comment. Math. Helv. {\bf 64}(1989), 412-461.
}

\newcommand{\Sififo}{\bibitem{Sififo} C.L. Siegel,
{\em \"{U}ber die Existenz einer Normalform analytischer Hamiltonscher
Differntialgleichungen in der N\"{a}he einer Gleichgewichtsl\"{o}sung},
Math. Ann. {\bf 128}(1954), 144-170.
}

\newcommand{\Brseon}{\bibitem{Brseon} A.D. Brjuno,
{\em Analytic form of differential equations I},
Trans. Moscow Math. Soc.
{\bf 25}(1971); {\em II}, 131--288 (1973).
}

\newcommand{\Stzeon}{\bibitem{Stzeon}
L. Stolovitch, {\em Singular complete integrability},
  Inst. Hautes \'Etudes Sci. Publ. Math.  No. 91 (2000),
133--210 (2001). }

\newcommand{\Mofifi}{\bibitem{Mofifi} J. Moser,
{\em
Nonexistence of integrals for canonical systems of differential equations},
Comm. Pure
Appl. Math. {\bf 8}(1955), 409--436.
}

\newcommand{\Mofiei}{\bibitem{Mofiei} J. Moser,
{\em On the generalization of a theorem of A. Liapounoff},
Comm. Pure Appl. Math.  {\bf 11}(1958), 257--271.
}

\begin{document}
%\begin{abstract}
%\end{abstract}
\maketitle

\setcounter{thm}{0}\setcounter{equation}{0}
\section{Introduction}

Let $\rr^{2n}$ be the standard symplectic space, equipped with
the symplectic two-form $\omega=dx_1\wedge dy_1+\dots+dx_n\wedge dy_n$.
Let $h(x,y)=O(2)$ be a real analytic functions defined near $0\in\rr^{2n}$.
Under suitable non-degeneracy condition on the quadratic form of $h$,
 one may find linear symplectic coordinates so that
\be{hxy=}
h(x,y)=\frac{1}{2}\sum_{1\leq j\leq\kappa}i\lambda_j(x_j^2+y_j^2)+\frac{1}{2}\Re
\sum_{\kappa<k, l\leq n,\lambda_l=\ov\lambda_k} \lambda_k
(x_k+ix_l)(y_k-iy_l)+O(3),
\end{equation}
where $\kappa=0,\dots, n$, and $\lambda_j$ is pure imaginary precisely when
$1\leq j\leq\kappa$, and $\lambda_1$, $-\lambda_1$, \dots, $\lambda_n$, $-\lambda_n$
are eigenvalues of $H_{zz}(0)J$ with $z=(x,y)$ and $Jx_j=y_j=-J^2y_j$.
One says that
$\lambda_1,\dots, \lambda_n$ are {\em non-resonant}, if
$\lambda\cdot\alpha\equiv\lambda_1\alpha_1+\dots+\lambda_n\alpha_n\neq0$
for all multi-indices of integers $\alpha\neq0$.
 The Birkhoff normal form says that under the non-resonance condition on $\lambda$,
 there is a formal symplectic transformation of $\rr^{2n}$
  sending $h$ into $\hat h$
 that is a real formal power series in $x_j^2+y_j^2$ $(1\leq j\leq \kappa)$,
 $(x_k+ix_l)(y_k-iy_l)$ ($\kappa<k,l\leq n$).
 Notice that, up to the order of
$\lambda_1,\dots, \lambda_n,
-\lambda_1,\dots, -\lambda_n$,
the Birkhoff normal form $\hat h$ is independent of the choice of the normalizing
 transformations.
 In~\cite{Sifoon}, Siegel showed that the Birkhoff normal form cannot be realized by convergent symplectic
 transformations in general. In fact, Siegel~\cite{Sififo} showed that
 when $\kappa=n\geq2$, for a real analytic
 function with any
 prescribed nonresonant $\lambda_1,\dots, \lambda_n$ and with generic
 higher order terms, there exists no convergent
 normalizing transformation.

 Despite  Siegel's divergence results and many other results,
 a basic question, which remains unsettled until now, is if
 there exists a divergent
 Birkhoff normal form arising from a real analytic function.
 This question was pointed out by Eliasson~\cite{Elnize}.
  To the author's knowledge, there seems  no example of divergent normal form
 in other normal form problems in the literature.
 The divergence of Birkhoff normal form implies, of course,
  that of all normalizing transformations
 of the given function. The importance of such a divergent normal form
 was demonstrated by P\'erez-Marco~\cite{Pezeth} very recently.

  In this paper we shall prove
  \bt{1} Let $\kappa=0,\dots, n$ and $n\geq2$. Assume
  that $\kappa\neq1$ when $n=2$.
  There exists a divergent Birkhoff normal
  form of some analytic real function \rea{hxy=}, defined near $0\in\rr^{2n}$
  and having non-resonant $\lambda_1,\dots,\lambda_n$.
\end{thm}
It is necessary to exclude the case
of non-real $\lambda_2/\lambda_1$ in the theorem when $n=2$. Indeed,
by a theorem of  Moser~\cite{Mofiei}, the Birkhoff normal form is always realized by some
   convergent
transformation when $n=2$ and $\lambda_2/{\lambda_1}$ is not real.
One can see, from the proof of the theorem, that the set of real analytic Hamiltonian functions
with divergent Birkhoff normal form is dense in a suitable topology. One may also apply
 a result of P\'erez-Marco~\cite{Pezeth} and the above theorem to conclude that
generic Hamiltonian functions
with the above quadratic form have divergent normal form too.

For the Birkhoff normal form theory, the reader is referred to, besides the above mentioned references,
 papers of Moser~\cite{Mofifi},
 R\"ussmann~\cite{Rusifo},~\cite{Rusise}, Brjuno~\cite{Brseon}, Vey~\cite{Veseei},
 Ito~\cite{Iteini}, Stolovitch~\cite{Stzeon},
 Giorgilli~\cite{Gizeon},
 and the author~\cite{zezemrl}, ~\cite{zeonaif}.
  Papers by  Brjuno~\cite{Brseon} and by P\'erez-Marco~\cite{Pezeth}
  contain extensive references also.

The proof of \rt{1} is based on the method of small divisors.
One would expect that the present
approach will have implications for other small-divisor problems.
We will however focus on the Hamiltonian functions, to demonstrate
how the small-divisors enter the
normal form.

\setcounter{thm}{0}\setcounter{equation}{0}
\section{Proof of the theorem}
We may restrict ourselves to $n=2$, since the sought $h$ for higher dimension can be
obtained trivially by adding suitable quadratic terms.

Consider a real analytic (real-valued) function
$$
h(x,y)=\sum_{j=1}^2\lambda_jx_jy_j+O(3),
$$
where $\lambda_1,\lambda_2$ are non-resonant.
Let $S(x,\hat y)$ be a real analytic function defined near $0\in\rr^2\times\rr^2$
with $S(x,\hat y)=O(d)$, $d>2$. Let $\varphi\colon
(x,y)\to (\hat x,\hat y)$ be a symplectic
map defined by
\be{hxj=}
\hat x_j=x_j-S_{\hat y_j}(x,\hat y), \quad \hat y_j=y_j+S_{x_j}(x,\hat y), \quad j=1,2.
\end{equation}
Note that
$$\varphi\colon\hat x_j=x_j-S_{y_j}(x,y)+O(d),
\quad \hat y_j=y_j+S_{x_j}(x,y)+O(d).
$$
Put $\hat h=h\circ \varphi^{-1}$.
Then $h(x,y)=\hat h(\hat x,\hat y)$ has the expansion
$$
\hat h(x,y)+\sum \lambda_j(x_jS_{x_j}(x,y)-y_jS_{y_j}(x,y))+O(d+1).
$$
Define the projection
$$\mathcal N\sum_{\alpha\beta} h_{\alpha\beta}x^\alpha y^\beta=\sum_{\alpha} h_{\alpha\alpha}x^\alpha y^\alpha.$$
Note that $h$ is in a Birkhoff normal form, if and only if $\mathcal Nh$
agrees with $h$.
For the special case of
 $h=\mathcal Nh+O(d)$ with $d\geq3$, taking
$$
S_{\alpha\beta}=\frac{1}{\lambda\cdot(\alpha-\beta)}h_{\alpha\beta}, \quad |\alpha|+|\beta|=d,
\quad\alpha\neq\beta
$$
yields
$\hat h=\mathcal N\hat h+O(d+1)$. In the above and in what follows
$\alpha,\beta$ stand for multi-indices of non-negative integers. We also write
$|\alpha-\beta|=|\alpha_1-\beta_1|+|\alpha_2-\beta_2|$.
 In general, inductively one finds
\be{salp}
S_{\alpha\beta}=
\frac{1}{\lambda\cdot(\alpha-\beta)}\{h_{\alpha\beta}+Q_{\alpha\beta}(h)\},
\quad \alpha\neq\beta,
\end{equation}
so that $\hat h_{\alpha\beta}=0$ for $\alpha\neq\beta$, i.e., so
that $\varphi$,  a formal
symplectic map of $\rr^4$, transforms $h$ into the Birkhoff normal form $\hat h$.
Notice that  the above expression $Q_{\alpha\beta}(h)$ stands for a
polynomial (with integer coefficients)
in quantities
$$
h_{\alpha'\beta'},
\frac{1}{\lambda\cdot(\alpha''-\beta'')};
\alpha''\neq\beta'',
\max\{|\alpha|'+|\beta'|,
 |\alpha''|+|\beta''|\}<|\alpha|+|\beta|.
$$
Note that $\lambda_1=h_{1,0,1,0}, \lambda_2=h_{0,1,0,1}$.
One also has
\be{halp}
\hat h_{\alpha\alpha}=
h_{\alpha\alpha}+D_{\alpha\alpha}(h),
\end{equation}
where $D_{\alpha\alpha}(h)$ is
a
polynomial
in quantities
$$
h_{\alpha'\beta'},
\frac{1}{\lambda\cdot(\alpha''-\beta'')}; \quad \alpha''\neq\beta'',
\max\{|\alpha|'+|\beta'|,
 |\alpha''|+|\beta''|\}<|\alpha|+|\beta|.
$$

We need to know more about the term $D_{\alpha\alpha}$ in \re{halp}.
\bl{2}
Let $S(x,\hat y)$ be a power series starting with terms of order $d$, and let $T=[S]_d$
be the sum of all monomials in $S$ of order $d>2$. Let $\varphi$ be the mapping defined
by \rea{hxj=}.
Let $\hat h=h\circ\varphi^{-1}$.
Assume that $h_{\alpha\beta}=0$ for
$|\alpha|+|\beta|<d$ and $\alpha\neq\beta$, and that $\hat h_{\alpha\beta}=0$
for $|\alpha|+|\beta|<2d-1$ and $\alpha\neq\beta$. Then
\be{dH}
\hat h(x,y)-\mathcal Nh(x,y)=\mathcal N\{\sum_{j,k=1}^2\lambda_j
(\frac{1}{2}T_{x_j}T_{y_j}+y_jT_{y_jy_k}T_{x_k}-x_jT_{x_jy_k}T_{x_k})
\}+O(2d-1).
\end{equation}
\el
\begin{proof}
Returning to \re{hxj=}, we get
\begin{align*}
   &\hat x_j=x_j-S_{y_j}(x,y)-\sum_{k=1}^2T_{y_jy_k}(x,y)T_{x_k}(x,y)+O(2d-2),  \\
   &\hat y_j=y_j+S_{x_j}(x,y)+\sum_{k=1}^2T_{x_jy_k}(x,y)T_{x_k}(x,y)+O(2d-2).
\end{align*}
Now
\begin{align*}
h(x,y)&=\hat h(\hat x,\hat y)=\sum_{j,k=1}^2\lambda_j(x_jT_{x_jy_k}T_{x_k}-y_jT_{y_jy_k}T_{x_k}
-\frac{1}{2}T_{x_j}T_{y_j})
\\
&+
\hat h(x,y)+\sum \alpha_j\hat h_{\alpha\alpha} x^{\alpha-\delta_j}y^{\alpha
-\delta_j}(x_jS_{x_j}-y_jS_{y_j})+O(2d-1),
\end{align*}
where $\delta_j=(0,\dots, 1,\dots, 0)$ with the $1$ at the $j$-th place.
Applying the projection $\mathcal N$ to both sides yields \re{dH}.
\end{proof}

The term $\lambda\cdot(\alpha-\beta)$ in \re{salp}
is the  small-divisor
 used by Siegel in his first proof~\cite{Sifoon}
  for the divergence of Birkhoff's normalization for
Hamiltonian
functions. Notably, this small divisor, when $|\alpha-\beta|=|\alpha|+|\beta|$,
 does not appear in \re{halp}.
We now identify the small-divisor that contributes to the
 divergence of a Birkhoff normal form.
\bl{3}
Keep nations and assumptions in \rl{2}. Let $N+m=d$, $\alpha=(N, m-1)$, $a=(N,0)$ and
$b=(0,m)$. Assume that $m\geq1$.
Then
\begin{align}\label{smal}
\hat h_{\alpha\alpha}&=h_{\alpha\alpha}-
\frac{m^2(\lambda_1N-\lambda_2)(h_{ab}+Q_{ab}(h))
(h_{ba}+Q_{ba}(h))}{(\lambda\cdot(a-b))^2}
\\&\quad
+\frac{ h_{ab}A_{N+m}(h)+h_{ba}B_{N+m}(h)}{\lambda\cdot(a-b)}
+C_{N+m}(h),\nonumber
\end{align}
where $A_{N+m}, B_{N+m}$ and $C_{N+m}$ are polynomials in
$
h_{\alpha'\beta'},
\frac{1}{\lambda\cdot(\alpha''-\beta'')}$
with
$$ \alpha''\neq\beta'', (\alpha'',\beta'')\neq(a,b), (b,a);
|\alpha'|+|\beta'|\leq N+m,
 |\alpha''|+|\beta''|\leq N+m.
$$
\el
\begin{proof}
Write
$$
T(x,y)=T_{ab}x_1^Ny_2^m+T_{ba}x_2^my_1^N+\sum_{(a',b')\neq (a,b), (b,a)}T_{a'b'}x^{a'}y^{b'}.
$$
Then we obtain
\begin{align*}
 \sum_{j,k}\lambda_jx_jT_{x_jy_k}T_{x_k}&=T_{ab}T_{ba}(
 \lambda_1Nm^2(x_1y_1)^{N}(x_2y_2)^{m-1}+\lambda_2mN^2(x_1y_1)^{N-1}(x_2y_2)^{m})+\dots,
   \\
 \sum_{j,k}\lambda_jy_jT_{y_jy_k}T_{x_k}&=0+\dots,\\
  \sum_j\lambda_jT_{x_j}T_{y_j}&=
T_{ab}T_{ba}(\lambda_1N^2(x_1y_1)^{N-1}(x_2y_2)^m+\lambda_2m^2(x_1y_1)^N(x_2y_2)^{m-1})
+\dots,
\end{align*}
where the omitted terms have coefficients that
are linear combinations with integer coefficients
 in $T_{ab}T_{a'b'}$,  $T_{ba}T_{a''b''}$, and
$T_{a'b'}T_{a''b''}$
with $(a',b')$, $(a'',b'')\neq(a,b)$, $(b,a)$.
Thus
\begin{align*}
 \hat h(x,y)-h(x,y)
&=
T_{ab}T_{ba}\{(\lambda_1-\lambda_2m)N^2(x_1y_1)^{N-1}(x_2y_2)^m\\
&\quad +(\lambda_2-\lambda_1N)m^2(x_1y_1)^N(x_2y_2)^{m-1}\}
+\dots.
\end{align*}
Combining \re{salp}, we obtain \re{smal}.
\end{proof}

\bp{4}
Let $h(x,y)=\lambda_1 x_1y_1+\lambda_2x_2y_2+O(3)$ be a real analytic function with $0<\lambda_1<\lambda_2$. Assume
that $\lambda_1,\lambda_2$ are non-resonant. Let $\varphi$ be any formal
symplectic map  so that $\hat h(x,y)=h\circ\varphi^{-1}(x,y)$
is in the Birkhoff normal form with
quadratic form $\lambda_1x_1y_1+\lambda_2x_2y_2$.
Then for $\alpha=(N,m-1)$, $a=(N,0), b=(0,m)$
with $m\geq1$, one has
\begin{align}\label{smal+}
\hat h_{\alpha\alpha}&=h_{\alpha\alpha}-
\frac{m^2(\lambda_1N-\lambda_2)(h_{ab}+Q_{ab}(h))
(h_{ba}+Q_{ba}(h))}{(\lambda\cdot(a-b))^2}
\\&\quad
+\frac{h_{ab}A_{ab}(h)+h_{ba}B_{ab}(h)+C_{ab}(h)}{\lambda\cdot(a-b)}
+\hat Q_{ab}(h),\nonumber
\end{align}
where  $Q_{ab}$ is a polynomial in
$
h_{\alpha'\beta'},
\frac{1}{\lambda\cdot(\alpha''-\beta'')}$
with
$$ \alpha''\neq\beta'', \max\{|\alpha'|+|\beta'|,
 |\alpha''|+|\beta''|\}<|a|+|b|,
$$
$\hat Q_{ab}$ is a polynomial in
$
h_{\alpha'\beta'},
\frac{1}{\lambda\cdot(\alpha''-\beta'')}$
with
$$ \alpha''\neq\beta'', (\alpha'',\beta'')\neq(a,b),(b,a),|\alpha'|+|\beta'|<2|\alpha|,
|\alpha''|+|\beta''|\leq |a|+|b|
$$
and  $ A_{ab}$, $B_{ab}, C_{ab}$  are polynomials in
$
h_{\alpha'\beta'}$, $
\frac{1}{\lambda\cdot(\alpha''-\beta'')}$
with
$$ \alpha''\neq\beta'', (\alpha'',\beta'')\neq(a,b),(b,a),
\max\{|\alpha'|+|\beta'|,
 |\alpha''|+|\beta''|\}\leq|a|+|b|.
$$
\ep
\begin{proof}
We apply a symplectic map $\varphi_1$ of the form \re{hxj=}, in which
$$
S(x,\hat y)=\sum_{\alpha\neq\beta,3\leq |\alpha|+|\beta|<N+m}
S_{\alpha\beta}x^\alpha\hat y^\beta,
$$
so that $\tilde h=h\circ\varphi_1^{-1}$ satisfies $\tilde h_{\alpha\beta}=0$
for all $\alpha\neq\beta$ and $|\alpha|+|\beta|<N+m$.
We know that
$\tilde h_{\alpha\beta}=h_{\alpha\beta}+D_{\alpha\beta}(h)$,
where $D_{\alpha\beta}(h)$ depends on $h_{\alpha'\beta'}$ with
$|\alpha'|+|\beta'|<|\alpha|+|\beta|$ and on $
1/{(\lambda\cdot(\alpha''-\beta''))}$ with $|\alpha''|+|\beta''|<|a|+|b|$,
$\alpha''\neq\beta''$.
Apply a formal symplectic map $\varphi_2$ of the form \re{hxj=} with
$$
S(x,\hat y)=\sum_{\alpha\neq\beta, |\alpha|+|\beta|\geq N+m}
S_{\alpha\beta}x^\alpha\hat y^\beta,
$$
so that $\tilde h\circ\varphi_2^{-1}$ is in the Birkhoff normal form.
By \re{smal}, in which $h$ is actually $\tilde h$ now, we can write (with abuse of
notation for $Q_{ab}(h)$)
\begin{gather*}
\tilde h_{\alpha\alpha}+C_{N+m}(\tilde h)=h_{\alpha\alpha}+\hat Q_{ab}(h),\\
\tilde h_{ab}+Q_{ab}(\tilde h)=h_{ab}+Q_{ab}(h),\\
\tilde h_{ab}A_{N+m}(\tilde h)+\tilde h_{ba}B_{N+m}(\tilde h)=
h_{ab}A_{ab}(h)+h_{ba}B_{ab}(h)+C_{ab}(h)
\end{gather*}
for $C_{ab}(h)=D_{ab}(\tilde h)A_{N+m}(\tilde h)+D_{ba}(\tilde h)B_{N+m}(\tilde h)$,
$A_{ab}(h)=A_{N+m}(\tilde h)$,
$B_{ab}(h)=B_{N+m}(\tilde h)$.

We have obtain \re{smal+}, via the above normalizing map $\varphi_2\varphi_1$. On the
other hand the Birkhoff normal form $\hat h$, with the same quadratic form as $h$,
 is independent of the normalizing map.
In other words, the right-hand side of \re{smal+} is independent of $\varphi$. Since each $\hat h_{\alpha\alpha}$ is
 a polynomial with integer coefficients in variables $h_{\alpha'\beta'}$, $\frac{1}{\lambda\cdot(\alpha''-\beta'')}$,
 we conclude that
each term in \re{smal+}  depends only on $h$ and is a polynomial in the sought form.
\end{proof}

We now restrict ourselves to
$|h_{\alpha\beta}|\leq2$ for all $\alpha,\beta$. Then
we have
\be{maxq}
\max\{|Q_{ab}|,|Q_{ba}|,|A_{ab}|,|B_{ab}|,|C_{ab}|,\hat Q_{ab}|\}\leq
\delta_{ab}(\lambda)^{-\tau_{ab}},
\end{equation}
where $\tau_{ab}>1$ is a constant independent of $\lambda$ and
$$
\delta_{ab}(\lambda)=\min\{\frac{1}{2},|\lambda\cdot(\alpha-\beta)|\colon\alpha\neq\beta,
|\alpha|+|\beta|\leq|a|+|b|, (\alpha,\beta)\neq(a,b),(b,a)\}.
$$
Put $\lambda_2=1$. Notice that for $a=(N,0), b=(0,m)$, one has $|a-b|=|a|+|b|$. Thus,
we can choose an irrational $\lambda_1\in(0,1)$ so that
\be{nmum}
|(a-b)\cdot\lambda|=|N\lambda_1-m|<\frac{\delta_{ab}(\lambda)^{\tau_{ab}}}{100(N+m)!},
\quad a=(N,0), b=(0,m)
\end{equation}
holds for a sequence $(N,m)=(N_j,m_j)$ with $N_j,m_j$ being positive integers. We may assume
that $N_{j+1}+m_{j+1}>2(N_j+m_j)$. Put
$a_j=(N_j,0), b_j=(0,m_j)$.

We now complete the proof of the theorem.

We construct $h$ for the case $\kappa=0$ first. We shall find $h$ whose coefficients $h_{\alpha\beta}$
 are real
and satisfy the extra condition $h_{\alpha\beta}=h_{\beta\alpha}$.
Put $h_{\alpha\beta}=0$
for all $\alpha,\beta$ with $|\alpha|+|\beta|>2$ and
$(\alpha,\beta)\neq(a_j,b_j),(b_j,a_j)$.
Inductively, we shall choose $h_{a_jb_j}=h_{b_ja_j}=0, 2$, or $-2$ as follows.
Notice that if $u_0,v_0$ are real and
 $|u_0v_0|<1$, then either $(u_0+2)(v_0+2)\geq2$ or $(u_0-2)(v_0-2)\geq2$;
 otherwise, we would have both  $u_0+v_0<-1/2$ and
 $u_0+v_0>1/2$, which is a contradiction.
Therefore for two real numbers $u_0,v_0$, choosing
 $(u,v)$ among $(0,0)$, $(2,2)$ and $(-2,-2)$ yields $|(u_0+u)(v_0+v)|\geq 1$.
This shows that  we can find $h_{a_jb_j}=h_{b_ja_j}=0$, $2$ or $-2$, so that
\be{habq}
 |(h_{ab}+Q_{ab}(h))(h_{ba}+Q_{ba}(h))|\geq 1, \quad a=a_j,b=b_j.
\end{equation}
Here, we already used  $N_{j+1}+m_{j+1}>2(N_j+m_j)$,
which implies that if \re{habq} holds for $a=a_j,b=b_j$
then it remains true no matter how $a_{j+1}, b_{j+1}$ are chosen. Now
 \re{smal+}-\re{habq} imply that for
$(N,m)=(N_j,m_j)$ and $|\lambda_1 N_j-1|>1$ we have
$$
|\hat h_{\alpha\alpha}|>\frac{m^2|\lambda_1 N-1|}{2|\lambda\cdot(a-b)|^2}>(N+m)!,
\quad \alpha=(N,m-1).
$$
This shows the divergence of $\hat h$.

We now construct $h$ for the case $\kappa=2=n$,
 via restricting the complexification of $h$ to a suitable totally
real subspace of $\cc^4$.

For the above analytic real function $h(x,y)$ on $\rr^2\times\rr^2$,
 its complexification,  denoted by  $h(z,w)$, is  holomorphic near $0\in\cc^2
\times \cc^2$.
Let $\varphi$ be a formal symplectic map of $\rr^4$, which is tangent to the identity,
so that $h\circ\varphi^{-1}(x,y)=g(x_1y_1,x_2y_2)$
is in the normal form.  Since $\varphi$ preserves $\omega=dx_1\wedge dy_1
+dx_2\wedge dy_2$,
its complexification,  still denoted by $\varphi$,
 preserves
$\omega^c=dz_1\wedge dw_1+dz_2\wedge dw_2$.

Let $L(\xi,\eta)=(\xi+i\eta,\xi-i\eta)$. Notice that $L^*\omega^c=-2i(d\xi_1\wedge d\eta_1+d\xi_2\wedge d\eta_2)$.
Thus
$\psi=L^{-1}\varphi L$ preserves $d\xi_1\wedge d\eta_1+d\xi_2\wedge d\eta_2$. Also
$\tilde h\circ\psi^{-1}(\xi,\eta)=g(\xi_1^2+\eta_1^2,\xi_2^2+\eta_2^2)$
for $\tilde h=h\circ L$. In other words,
$\tilde h\circ\psi^{-1}$ is a (formal holomorphic) Birkhoff normal form with respect to
the  holomorphic symplectic $2$-form $d\xi_1\wedge d\eta_1+d\xi_2\wedge d\eta_2$.
Notice that the quadratic form of $\tilde h$
 is now $\lambda_1(\xi_1^2+\eta_1^2)+\lambda_2(\xi_2^2+\eta_2^2)$.
Let $e$ be the restriction of $\tilde h$
on $\rr^2\times\rr^2\colon \xi=\ov\xi,\eta=\ov\eta$.
Since $h_{\alpha\beta}=\ov h_{\beta\alpha}$ by  construction, then $e$ is real-valued.
Thus $e(\xi,\eta)$ is an analytic real function
of the form $\lambda_1(\xi_1^2+\eta_1^2)+\lambda_2(\xi_2^2
+\eta_2^2)+O(3)$, while $L^*\omega^c$,
restricted to $\rr^2\times\rr^2\colon\xi=\ov\xi,\eta=\ov\eta$,
is a constant multiple of the standard symplectic real $2$-form.
Therefore $\tilde h\circ\psi^{-1}$, restricted to $\xi=\ov\xi,\eta=\ov\eta$, is a real
Birkhoff normal
form of $e$;
since $h$ diverges,  one readily sees the divergence of
the restriction.

\bibliographystyle{plain}

\end{document}